\documentclass[12pt]{article}
\usepackage{amsmath, amssymb}

\begin{document}

\title{
REALIZATION OF THE \\
RIEMANN HYPOTHESIS \\
VIA \\
COUPLING CONSTANT SPECTRUM}

\author{R. Acharya \footnote{I dedicate this note to my teacher, George Sudarshan
of the University of Texas at Austin}\\
Department of Physics\\
Arizona State University\\
Tempe, AZ 85287-1504}

\maketitle
\thispagestyle{empty}

\begin{abstract}
We present a Non-relativistic Quantum mechanical model, which exhibits the
realization of Riemann Conjecture. The technique depends on
exposing the $S$-wave Jost function at zero energy and in
identifying it with the Riemann $\xi(s)$ function following a
seminal paper of N. N. Khuri.
\end{abstract}

\newpage
\setcounter{page}{1}

We begin by recalling the all-too-familiar lore that the Riemann
hypothesis has been the Holy Grail of mathematics and physics for
more than a century \cite{Borwein}. It asserts that \textbf{all} the zeros of
$\xi(s)$ have $\sigma = \frac{1}{2}$, where $s=\sigma \pm it_{n}$,
$n = 1, 2, 3 \dots \infty$. It is believed all zeros of $\xi(s)$
are simple. The function $\zeta(s)$ is related to the Riemann
$\xi(s)$ function via the defining relation \cite{Borwein},

\begin{equation}\label{eq:1}
\xi(s) = \frac{1}{2} s(s-1)\pi^{-\frac{S}{2}} \Gamma \Bigl( \frac{s}{2} \Bigr)
\zeta(s)
\end{equation}
so that $\xi(s)$ is an entire function, where

\begin{equation}\label{eq:2}
\zeta(s) = \sum^{\infty}_{n=1} n^{-s}, \quad s = \sigma + it, \quad \sigma >
1
\end{equation}

$\zeta(s)$ is holomorphic for $\sigma > 1$ and can have No zeros for
$\sigma > 1$. Since $1/\Gamma(z)$ is entire, the function $\Gamma
\bigl( \frac{s}{2} \bigr)$ is non-vanishing, it is clear that
$\xi(s)$ \textbf{also} has no zeros in $\sigma > 1$: the zeros
of $\xi(s)$ are confined to the ``critical strip" $0 \leqslant \sigma
\leqslant 1$. Moreover, if $\rho$ is a zero of $\xi(s)$, then so is
$1-\rho$ and since $\overline{\xi(s)} = \xi(\overline{s})$, one deduces that
$\overline{\rho}$ and $1 - \overline{\rho}$ are also zeros. Thus the Riemann
zeros are symmetrically arranged about the real axis
\textbf{and} also about the ``critical line" given by $\sigma =
\frac{1}{2}$. The Riemann Hypothesis, then, asserts that ALL zeros
of $\xi(s)$ have Re $s= \sigma = \frac{1}{2}$.

We conclude this introductory, well-known remarks with the
assertion that every entire function $f(z)$ of \textbf{order
one} and ``\textbf{infinite type}" (which guarantees the
existence of \textbf{infinitely many Non-zero}
zeros can be represented by the Hadamard factorization, to wit \cite{Borwein},

\begin{equation}\label{eq:3}
f(z) = z^{m} e^{A} e^{BZ} \prod^{\infty}_{n=1}
\Bigl( 1-\frac{z}{z_{n}} \Bigr) \exp \Bigl( \frac{z}{z_{n}} \Bigr)
\end{equation}
where `m' is the multiplicity of the zeros (so that $m=0$, for
simple zero).

Finally, $\xi(s) = \xi(1-s)$ is indeed an entire function of order
one and infinite type and it has \textbf{No zeros either for
$\mathbf{\sigma > 1}$ or $\mathbf{\sigma < 0}$.}

We begin our brief note by defining the non-relativistic, quantum
mechanical potential model in 3 dimensions.

The zero energy (i.e., $k^{2} = 0$, $\frac{2m}{\hbar^{2}} = 1$)
Schrodinger equation reads [for \textbf{zero} angular momentum,
$l = 0$ ($S$-waves)]

\begin{equation}\label{eq:4}
\frac{1}{r^{2}} \frac{\partial}{\partial r} \biggl( r^{2}
\frac{\partial \psi(r)}{\partial r} \biggr) + \biggl(
0 - \frac{\lambda}{r^2} \biggr) \psi(r) = 0
\end{equation}

\textbf{The potential we choose is:}
\begin{equation*}
\tag{\ref{eq:4}$'$}
\begin{split}
V(r) = \frac{\lambda}{r^{2}}
\end{split}
\end{equation*}

The canonical change of the independent variable $\psi(r)$ to
$\frac{U(r)}{r}$ results in the equation, well-known to every one
(!):

\begin{equation}\label{eq:5}
U''(r) - \frac{\lambda}{r^{2}} U(r) = 0
\end{equation}

Eq. (\ref{eq:5}) has the solution:

\begin{equation}\label{eq:6}
U(r) = r^{s} \Rightarrow \psi(r) = r^{s-1}
\end{equation}
where $s$ is \textbf{constrained} to satisfy the relation,

\begin{equation}\label{eq:7}
\lambda = s (s-1)
\end{equation}

For the Regular solution at the origin, i.e., $\psi(0) = 0$, we
require that Real $s > 1$. Eq. (\ref{eq:7}) will play a crucial role in the
upcoming analysis. We note in passing that the Repulsive nature of
the real potential, i.e., $V(r) = \frac{\lambda}{r^{2}}$,
$\lambda$ real and $> 0$, requires that there are \textbf{No
bound states.} Further due to the scale invariance of the problem,
i.e., the potential is scale-free (and we are suppressing the
overall factor $\frac{1}{M}$ in $V(r)$). The \textbf{coupling
constant $\mathbf{\lambda}$ is dimensionless.}

The inverse-square potential acts like a centrifugal term in the
free Schro-dinger equation. It is a straightforward exercise in
undergraduate physics to determine the corresponding $S$-wave
phase shift for scattering solution \cite{Barford}:

\begin{equation}\label{eq:8}
\delta = \frac{\pi}{4} - \frac{\pi\nu}{2}
\end{equation}
where

\begin{equation}\label{eq:9}
\nu = \sqrt{\frac{1}{4} + \lambda} > \frac{1}{2}
\end{equation}

We note in passing that for the free case, $\lambda = 0$, $\nu =
\frac{1}{2}$ and the phase shift does indeed vanish and the $S$
matrix is unity:
\begin{equation}\label{eq:10}
S = \exp(2i\delta)
\end{equation}

It is important to emphasize that the phase shift $\delta$ in Eq.
(\ref{eq:8}) and hence also the $l = 0$ partial wave scattering
matrix in Eq. (\ref{eq:10}) are \textbf{independent of energy} and
the $S$ matrix is thus manifestly scale invariant. (There are no
bound states).

We go on to summarize the seminal idea of Khuri (and Chadan)
\cite{Khuri}. It is well known that the $S$ wave Jost function is
identified via the defining relation,

\begin{equation}\label{eq:11}
S(k) = \frac{F_{-}(k)}{F_{+}(k)}
\end{equation}
where $F_{+}(k)$ is holomorphic in the upper half complex $k$
plane and $F_{-}(k)$ is holomorphic in the lower half of $k$
plane, where $k$ is the momentum (in the center-of-mass frame).

Thus,

\begin{equation}\label{eq:12}
\begin{split}
F_{+}(k) = {S}^{-1} (k) F_{-}(k), \qquad  & Im k = 0 \\
                                         & -\infty < k < \infty
\end{split}
\end{equation}
and

\begin{equation}\label{eq:13}
\det S(k) \neq 0, \qquad Imk = 0
\end{equation}

Khuri's idea is the following: the $S$-wave Jost function for a
potential

\begin{equation}\label{eq:14}
\lambda V^{\ast} = \frac{\lambda}{r^{2}}
\end{equation}
is an entire function of $\lambda$ with an infinite number of
zeros extending to infinity. For a repulsive potential $V$ and at
\textbf{zero energy}, these zeros of the coupling constant
$\lambda$, will all be \textbf{real} and \textbf{negative}
(see the elaboration of this fact later on), i.e.,
$\lambda_{n}(k^{2} = 0) < 0$.

By rescaling $\lambda$ such that $\lambda_{n} < -\frac{1}{4}$ and
changing variables to $s$, with $\lambda = s(s-1)$, it follows
that as a function of $s$, the $S$-wave Jost function has
\textbf{only} (!) zeros on the line $s_{n} = \frac{1}{2} +
it_{n}$. Thus, if we can ``find" a \textbf{repulsive} potential
$V^{\ast}$ whose coupling constant spectrum coincides with the Riemann
zeros of $\xi(s)$, this will unambiguously establish the holy
grail of physics and mathematics! Khuri commented that
\textbf{``This will be a very difficult and unguided search."}
That said, we wish to make the crucial observation: For the
potential at hand, i.e., $V(r) = \lambda V^{\ast}(r)$,
$V^{\ast}(r)=\frac{1}{r^{2}}$, $\lambda$ real \& positive, the
coupling constant $\lambda$ must \textbf{necessarily} satisfy
the equality,

\begin{equation}\label{eq:15}
\lambda = s(s-1)
\end{equation}

In other words, there is ``no need to change variables!" We
proceed by summarizing the ``solution" to the boundary value
problem defined by Eq. (\ref{eq:12}) for the $S$ wave Jost function.

In the case of uncoupled partial waves the solution to the
boundary value problem formulated via Eq. (\ref{eq:12}) has been
given by Krutov, Muravyev and Troitsky \cite{Krutov} in 1996. The
``Solution" reads:

\begin{equation}\label{eq:16}
F_{\pm}(k) = \Pi_{\pm}(k) \exp \biggl(
\frac{1}{2\pi i} \int^{\infty}_{-\infty}
\frac{\ln(S^{-1}(k) \Pi^{2}_{-}(k))}{k'-k\mp io}dk'
\biggr)
\end{equation}
where

\begin{equation}\label{eq:17}
\Pi_{\pm}(k) = \prod^{m'}_{j=1} \frac{k \mp ik_{j}}{k \pm ik_{j}},
\quad k_{j} >0, \quad j=1,2,\dots,m' \quad (m' < \infty)
\end{equation}

\begin{equation}\label{eq:18}
\Pi_{\pm}(k) \equiv 1, \quad m' = 0
\end{equation}
where $m'$ is the number of bound states.

It is well-known, of course that

\begin{equation}\label{eq:19}
\frac{1}{k'-k \mp io} = P \frac{1}{k'-k} \pm i\pi \delta(k'-k)
\end{equation}

Proceeding further, we set the Jost function ($S$ wave), at zero
energy

\begin{equation}\label{eq:20}
F_{+}(s(s-1); k^{2} = 0 ) \equiv \chi(s)
\end{equation}

The simplification of Eq. (\ref{eq:16}) is immediate. One finds
easily that

\begin{equation}\label{eq:21}
F_{+} = \chi(s) = \exp(-2i\delta), \quad S(s) = \exp(2i\delta)
\quad \text{and} \quad F_{-}(s) = 1
\end{equation}

We note in passing that since the partial wave $S$ matrix is
independent of $k$ (energy) due to scale invariance, one can set
$k=0$ in simplifying Eq. (\ref{eq:16}). Recall that there are no
bound states either, $m' = 0$ and $\Pi_{\pm}(k) \equiv 1$.

Since $\chi(s)$ defined in Eq. (\ref{eq:20}) is entire \cite{Levin}
and has \textbf{only zeros} on the real line Re $s_{n} =
\frac{1}{2}$ (see clarification of this later) and the number of
non-zero zeros are infinite on the ``critical" line Re $s_{n} =
\frac{1}{2}$, the Hadamard factorization Eq. (\ref{eq:3}) is valid,
i.e.,

\begin{equation}\label{eq:22}
\chi(s) = e^{A} e^{Bs} \prod^{\infty}_{n=1} \biggl( 1-\frac{s}{s_{n}}
\biggr)\exp\biggl(\frac{s}{s_{n}}\biggr)
\end{equation}

It will turn out that both the constants $A$ and $B$ can be
identified (see below) and also that $\chi(s)$ has No zeros either
for $s > 1$ or $s < 0$.

From Eq. (\ref{eq:8}), Eq. (\ref{eq:9}) and Eq. (\ref{eq:10}), we
obtain

\begin{equation}\label{eq:23}
\nu = \frac{1}{2} - \frac{2}{\pi} \delta = \frac{1}{2} +
\frac{1}{\pi} |\ln\chi(s)| > \frac{1}{2}
\end{equation}
where we have set (Eq. (\ref{eq:21}))

\begin{equation}\label{eq:24}
\ln \chi(s) = +i|\ln \chi(s)|
\end{equation}
in order to satisfy that $\delta < 0$

\begin{equation}\label{eq:25}
\delta = - \frac{1}{2} |\ln \chi(s)|
\end{equation}

[Recall that for a \textbf{repulsive} potential, $\lambda > 0$ the
phase shift has to be \textbf{negative}!]

\begin{equation}\label{eq:26}
\therefore \lambda = s(s-1) = \nu^{2} - \frac{1}{4}
\end{equation}
giving [From Eq. (\ref{eq:23}), Eq. (\ref{eq:25})]

\begin{equation}\label{eq:27}
\lambda = \frac{1}{\pi} | \ln \chi(s) | + \frac{1}{\pi^{2}} |\ln
\chi(s)|^{2} > 0
\end{equation}
where $\lambda$ is real positive and $s$ real and $> 1$.

From Eq. (\ref{eq:26}) and Eq. (\ref{eq:27}), we obtain $\mathbf{\chi(s)}$
\textbf{for real $\mathbf{s > 1}$:}

\begin{equation}\label{eq:28}
|\ln \chi(s)| = \pi(s-1) > 0 \qquad s \text{, real \&} > 1
\end{equation}

\begin{align}
\therefore \chi(s) & = e^{\pi(s-1)}, \quad \mathbf{s > 1} \textbf{,
real} \label{eq:29} \\
                   & = e^{-\pi s}, \qquad \mathbf{s < 0} \textbf{,
                   real} \label{eq:30}
\end{align}

$A = - \pi$, $B = \pi$ for $s > 1$, real (from Eq. (\ref{eq:22}))

Notice, as promised that $\chi(s)$ is entire (and order one) and
\textbf{has No zeros for $\mathbf{s > 1}$ and $\mathbf{s < 0}$.}
It is worth noticing that since $\lambda = s(s-1)$ is invariant
under $s \rightarrow 1-s$, the forms of $\chi(s)$ in Eq.
(\ref{eq:29}) and Eq. ({\ref{eq:30}) are consistent.

We now comment on the \textbf{all-important} restriction that for
a real, repulsive potential the $S$-wave Jost function $\lambda V$
is an entire function of $\lambda$ with an infinite number of
zeros extending to infinity and these zeros will all be
\textbf{real and negative}. [This can be verified following Khuri
{See his Eq. (1.1) and Eq. (1.2)] i.e.,

\begin{equation}\label{eq:31}
[I m \lambda_{n} (i \tau)] \int^{\infty}_{0} | f (\lambda_{n}(i \tau));
i \tau; r |^{2} dr = 0
\end{equation}
where $k = i \tau$.

We observe that Eq. (\ref{eq:31}) does \textbf{not} have the
potential $V(r)$ under the integration sign. This comes about by
examining the Schrodinger equation for $f(k,r)$ [which is
identical to the case of one dimension] and its complex conjugate
$f^{*}(k,r)$, and by performing the ``canonical" operation, i.e.,
multiplying the equation for $f(k,r)$ by $f^{*}(k,r)$,
subtracting and then importantly multiplying the resulting
equation by $r^{2}$. This ``gets rid of" the potential term on
R.H.S. Subsequently, one integrates by parts the equation thus
obtained, making use of the fact that the Wronskian of $f(k,r)$
and $f^{*}(k,r)$ is \textbf{independent} of $r$ (because they
satisfy the \textbf{same} Schrodinger equation: See Chaden and
Sabatier, \cite{Khuri}). The end result is Eq. (\ref{eq:31})
above, the potential $V(r)$ \textbf{not} present (in contrast to
Khuri's Eq. (1.2)). The rationale behind this is to ensure that
the integral in Eq. (\ref{eq:31}) is \textbf{now finite}, under
the stated conditions (See Eq. (36$'$) and Eq. (36$''$). If the
potential $V(r) = \lambda \frac{1}{r^{2}}$ were present, then the
corresponding integral in Eq. (\ref{eq:31}) would diverge!

The crucial point to verify is that the integral in Eq.
(\ref{eq:31}) is \textbf{non vanishing} and \textbf{finite. Only
then,} can one conclude that

\begin{equation}\label{eq:32}
I m \lambda_{n} (i \tau) = 0.
\end{equation}

We proceed to demonstrate this as follows.

The $S$ wave Jost \textbf{solution} $f(k,r)$ for the
potential $V(r) = \frac{\lambda}{r^{2}}$ ($\lambda >
0$) is well-known \cite{Barford}:

\begin{equation}\label{eq:33}
f(k,r) = \sqrt{\frac{\pi k r}{2}} e^{i(\frac{\pi}{2}\nu +
\frac{\pi}{4})} H^{(1)}_{\nu} (kr)
\end{equation}
where $H^{(1)}_{\nu} (kr)$ is Hankel function of $I$ kind.

One \textbf{cannot} obtain the Jost function by taking the limit
$r \rightarrow 0$ in Eq. (\ref{eq:33}) because this limit does
\textbf{not} exist. This only works for REGULAR POTENTIALS! (See
Chadan and Sabatier, Page 10, Eq. (I.3.6).)
We can bypass this conundrum, because we
already ``know" the $S$ wave Jost function $F_{+}(k) = F_{+}(0)$
from Eqs. (\ref{eq:8}), (\ref{eq:9}) and (\ref{eq:21}).

We observe in passing that the Jost solution Eq. (\ref{eq:33}) has
the required asymptotic behavior, i.e.,

\begin{equation}\label{eq:34}
f(k,r) \mathrel{\mathop{\longrightarrow}\limits_{r\rightarrow \infty}} e^{ikr}
\end{equation}

Since
\begin{equation}\label{eq:35}
H^{(1)}_{\nu}(z) \sim \sqrt{\frac{2}{\pi z}}
e^{i(z-\frac{\pi}{2} \nu - \frac{\pi}{4})},
\end{equation}
it is now straightforward to check that
\begin{equation}\label{eq:36}
\int^{\infty}_{0} |f(\lambda_{n}(i\tau); i\tau; r|^{2} dr < \infty
\end{equation}
by plugging in Eq. (\ref{eq:33}) for $f(k,r)$.

\begin{equation*}
\tag*{Eq. (4$'$)}
\begin{split}
V(r) = \frac{\lambda}{r^{2}}
\end{split}
\end{equation*}

Eq. (\ref{eq:31}), (\ref{eq:33}) give

\begin{equation*}
[ I m \lambda_{n} (i\tau)\tau ] \frac{2}{\pi}
\int^{\infty}_{0} r K^{2}_{\nu} (\tau r) dr = 0
\tag{\ref{eq:36}$'$}
\end{equation*}
where
\begin{equation*}
\int^{\infty}_{0} r K^{2}_{\nu} (\tau r) dr = \frac{1}{8}
\frac{1}{\tau^{2}} \frac{\pi \nu}{\sin \pi \nu}
\tag{\ref{eq:36}$''$}
\end{equation*}

$\mathbf{\therefore}$ \textbf{Eq. (\ref{eq:32}) follows from Eq. (36$'$), for $\mathbf{\tau
\neq 0}$ ($\mathbf{\tau > 0}$) \cite{Gradshteyn}}

\begin{equation*}
I m \lambda_{n} (i \tau) = 0, \qquad \nu > \frac{1}{2}, \quad \nu \neq 1, 2, 3 \dots \infty. \tag{32}
\end{equation*}

We now follow Khuri \cite{Khuri}:

All zeros of $\lambda_{n}(i\tau)$ must \textbf{be real} ($\tau >
0$) ``But $\lambda_{n}(i\tau)$ must be negative since the
potential $[\lambda_{n}(i \tau) V^{\ast}]$ will have a bound state
at $E = - \tau^{2}$ and that could not happen if $V^{\ast}
\geqslant 0$ and $\lambda_{n}(i \tau) > 0$. Hence by continuity,
$\mathbf{\lambda_{n}(0)}$\textbf{, for all $\mathbf{n}$, is real and
negative.} The zero energy coupling spectrum lies on the negative real
line for $V^{\ast} \geqslant 0$." In our case, in view of Eq.
(\ref{eq:7}), i.e.,
\begin{center}
$\lambda = s(s-1)$ is mandatory!
\end{center}

\textbf{Thus, we conclude unequivocally that $\mathbf{\chi(s)}$
has all its infinite number of nonzero zeros on the critical
line. And can be set equal to Riemann's $\mathbf{\xi}$
function!}

This establishes the Holy grail of physics and mathematics!!

We conclude by commenting on the precise relation between the $S$
wave Jost function, $\chi(s)$ and Riemann's $\xi(s)$ function.
Both are entire, order one and infinite type and possess identical
infinite number of non-zero zeros \textbf{only} on the critical
line Re $s = \frac{1}{2}$ and both have no zeros either for
$\sigma > 1$ or $\sigma < 0$. Their ratio, by Hadamard's theorem
is given by

\begin{equation}\label{eq:37}
\frac{\chi (s)}{\xi(s)} = e^{\alpha + \beta s}, \quad \text{all} \; s
\end{equation}

Since \cite{Mezzadri}

\begin{equation}\label{eq:38}
\xi (s) = \frac{1}{2} e^{b_{0}s} \prod_{n} \biggl(
1-\frac{s}{s_{n}} \biggr) e^{\frac{s}{s_{n}}}, \quad \textbf{all}
\; s
\end{equation}
where

\begin{equation}\label{eq:39}
b_{0} = - \frac{1}{2} \gamma - 1 - \frac{1}{2} \ln 4 \pi
\end{equation}
\begin{center}
($\gamma$ is Euler's constant)
\end{center}

We obtain the identification,

\begin{align}
          \alpha & = \ln2 - \pi \label{eq:40} \\
\text{and} \quad \beta & = \pi + \frac{1}{2} \gamma + 1 - \frac{1}{2}
\ln 4\pi \label{eq:41}
\end{align}

I wish to thank Professor Coffey, Professor Odlyzko and Professor
Ivic for a communication on the zeros of $\xi(s)$. I am especially
grateful to Professor Okubo \cite{Okubo} for informing about his
work on the Riemann Hypothesis. I am grateful for a communication
from Professor Carlos Castro and Professor H. C. Rosu. I am grateful to Irina Long for
her generous help in the preparation of this paper.

\end{document}